\documentclass[12pt]{amsart}
\usepackage{amsmath}
\usepackage{amssymb}
\usepackage{stmaryrd}
\usepackage{epsfig,color,esint}
\usepackage{blindtext}
\usepackage{enumerate}
\usepackage{enumitem}
\usepackage{hyperref}
\usepackage{graphicx}
\usepackage{url}
\usepackage{bbm}
\usepackage{nicefrac,mathtools}
\usepackage{mathrsfs}
\usepackage{cite}
\usepackage{comment}
\usepackage{pgfplots}
\usepackage{tikz}
\usetikzlibrary{arrows,shapes,trees,backgrounds}

\numberwithin{equation}{section}

\newtheorem{prop}{Proposition}[section]
\newtheorem{theo}{Theorem}[section]
\newtheorem{lemm}{Lemma}[section]
\newtheorem{coro}{Corollary}[section]
\newtheorem{rema}{Remark}[section]

\newtheorem{defi}{Definition}[section]

\def\begeq{\begin{equation}}
\def\endeq{\end{equation}}

\begin{document}

\title{Nonuniqueness and Nonexistence Results for the $L_p$-dual Minkowski Problem with Supercritical Exponents}
\author{Shi-Zhong Du}
\author{YanNan Liu}
\author{Jian Lu}
\thanks{The author Du was partially supported by Natural Science Foundation of China (12171299).}
\thanks{The author Liu was partially supported by Natural Science Foundation of China (12071017 and 12141103), and Natural Science Foundation of Beijing Municipality (1212002).}
\thanks{The author Lu was partially supported by Natural Science Foundation of China (12122106).}
  \address{The Department of Mathematics,
            Shantou University, Shantou, 515063, P. R. China.} \email{szdu@stu.edu.cn}
  \address{School of Mathematics and Statistics,
            Beijing Technology and Business University, Beijing, 100048, P. R. China.} \email{liuyn@th.btbu.edu.cn}
  \address{School of Mathematics Sciences,
            South China Normal University, Guangzhou, 510631, P. R. China.} \email{lj-tshu04@163.com}

\renewcommand{\subjclassname}{%
  \textup{2010} Mathematics Subject Classification}
\date{Feb. 2024}

\begin{abstract}
   In this paper, the $L_p$-dual Minkowski problem
     $$
       \det(\nabla^2h+hI)=(h^2+|\nabla h|^2)^{\frac{n-q}{2}}h^{p-1}f, \ \ \forall x\in{\mathbb{S}}^{n-1}
     $$
of Monge-Amp\`{e}re type were studied for different $p$ and $q$. Some new nonuniqueness results were obtained for the range $p\leq-n+1, p<q-\lambda_1(n,k)$ and $f\equiv1$, where $\lambda_1(n,k)$ is the best constant of the Poincar\'{e} inequality on ${\mathbb{S}}^{n-1}$ with $k-$symmetricity. The second part of this paper is devoted to prove some new nonexistence results for the supercritical range $p\leq-q, q\geq n$ on all dimensional spaces. The key ingredient of our proof was based on a generalization of Chou-Wang type identity obtained in \cite{chou2006lp} for $q=n, p=-q$ to a full range of $(p,q)$.
\end{abstract}

\maketitle\markboth{$L_p$-dual Minkowski problem}{$k-$symmetricity}

\tableofcontents

\section{Introduction}
The Minkowski problem is to determine a convex body with prescribed curvature or other similar geometric data. It plays a central role in the theory of convex bodies. Various Minkowski problems \cite{aleksandrov1942existence,chen2018planar,cheng1976regularity,guan2012hypersurfaces,guan2003christoffel,huang2016geometric,pogorelov1973extrinsic,guang2022l_p} have been studied extensively after Lutwak \cite{lutwak1975dual,lutwak1993brunn}, who proposes two variants of the Brunn-Minkowski theory including the dual Brunn-Minkowski theory and the  $L_p$ Brunn-Minkowski theory. Besides, there are singular cases such as the logarithmic Minkowski problem and the centro-affine Minkowski problem \cite{boroczky2012log,chou2006lp}.

Given a convex body $\Omega\in{\mathbb{R}}^{n}$ containing the origin, we denote by $\Omega\in{\mathcal{K}}_0$. The support function $h: {\mathbb{S}}^{n-1}\to{\mathbb{R}}$ of $\Omega$ is defined by $h(x)\equiv\max\big\{z\cdot x\big|\ z\in\Omega\big\}$ for each $x\in{\mathbb{S}}^n$. For regular convex body $\Omega$, $[A_{ij}]\equiv[h_{ij}+h\delta_{ij}]$ is a positive definite matrix, where $h_{ij}=\nabla^2_{ij}h$ stands for the Hessian matrix of $h$ acting on an orthonormal frame $\{e_i\}_{i=1}^n$ of ${\mathbb{S}}^{n-1}$ and $\delta_{ij}$ denotes the Kronecker delta symbol. As usually, we also use the notation $[A^{ij}]$ to denote the inverse matrix of $[A_{ij}]$, and the notation $[U^{ij}]$ to denote the cofactor matrix of $[A_{ij}]$. Letting
   $$
    r(y)\equiv\sup\big\{\lambda>0|\ \lambda y\in\Omega\big\}, \ \forall y\in{\mathbb{R}}^n
   $$
be the radial boundary distance function of $\Omega$, if one sets $Z(y)\equiv r(y)y\in\partial\Omega, y\in{\mathbb{S}}^{n-1}$, there holds
    $$
     Z=h(x)x+\nabla h(x), \ \ \forall x\in{\mathbb{S}}^{n-1}.
    $$
Let us introduce the spherical image mapping of boundary $\partial\Omega$ by $\nu(Z)\equiv\big\{x\in{\mathbb{S}}^{n-1}\big|\ z\cdot x=h(x)\big\}$. Then, one can define the radial Gauss mapping by
   $$
    G(\omega)\equiv\big\{\nu(Z(y))\big|\ y\in\omega\subset{\mathbb{S}}^{n-1}\big\}
   $$
and the reverse radial Gauss mapping by
   $$
    G^{-1}(\omega)\equiv\big\{y\in{\mathbb{S}}^{n-1}\big|\ \nu(Z(y))\in\omega\big\}.
   $$
For convex body $\Omega\in{\mathcal{K}}_0$ and two parameters $p,q\in{\mathbb{R}}$, the $L_p$ dual curvature measure introduced by Lutwak-Yang-Zhang \cite{lutwak2002sharp} is defined by
  $$
   C_{p,q}(\Omega,\omega)\equiv\int_{G^{-1}(\omega)}\frac{r^q(y)}{u^p(G(y))}dy.
  $$
As proposed in \cite{lutwak2002sharp}, given a finite Borel measure $\mu$ on ${\mathbb{S}}^{n-1}$, the $L_p$ dual Minkowski problem is to find a convex body $\Omega\in{\mathcal{K}}_0$ such that
   \begin{equation}\label{e1.1}
    dC_{p,q}(\Omega,\cdot)=d\mu.
   \end{equation}
On the smooth category, solving the $L_p$ dual Minkowski problem \eqref{e1.1} is boiled down to solve the fully nonlinear equation
  \begin{equation}\label{e1.2}
    \det(\nabla^2h+hI)=fh^{p-1}(h^2+|\nabla h|^2)^{\frac{n-q}{2}}, \ \ \forall x\in{\mathbb{S}}^{n-1}
  \end{equation}
of Monge-Amp\`{e}re type. The dual Minkowski problem \eqref{e1.2} contains the $L_p$ Aleksandrov problem proposed by Huang-Lutwak-Yang-Zhang in \cite{huang2018lp}, the logarithmic Minkowski problem \cite{boroczky2013logarithmic,zhu2014logarithmic}, the centro-affine Minkowski problem \cite{chou2006lp,zhu2015lp} and other unsolved problems. When $f\equiv1$, \eqref{e1.2} was reduced to
    \begin{equation}\label{e1.3}
     \det(\nabla^2h+hI)=(h^2+|\nabla h|^2)^{\frac{n-q}{2}}h^{p-1}, \ \ \forall x\in{\mathbb{S}}^{n-1}.
    \end{equation}
The issue of uniqueness of solution of \eqref{e1.2} or \eqref{e1.3} has attracted much attentions. As one knows, the uniqueness of \eqref{e1.2} was shown in \cite{huang2018lp} for $q<p$ and was shown in \cite{chen2021lp} for $q=p$ up to a constant. In the case $-2\leq p<q\leq\min\{2, p+2\}$, Chen-Huang-Zhao \cite{chen2019smooth} obtained the uniqueness of even solution of \eqref{e1.3}. When $1<p<q\leq2$, Chen-Li \cite{chen2021lp} showed uniqueness without evenness. On another direction, the nonuniqueness of \eqref{e1.3} on the planar case has been obtained by Huang-Jiang \cite{huang2019variational} for $p=0$ and even $q\geq6$, and by Chen-Chen-Li \cite{chen2021variations} for $pq\geq0, q-p>4$ or $p<0<q, 1+\frac{1}{p}<\frac{1}{q}$ and $q-p>4$. On the higher dimensional case, Li-Liu-Lu \cite{Li2022nonuniqueness} have shown nonuniqueness for $p<0<q$ and additional restriction. In \cite{yongsheng2021multiple}, Jiang-Wang-Wu showed a result of multiple solutions for $p<0$ and even $q\geq2$ such that $q-p>16$. More recently, the counting number of the solutions of \eqref{e1.3} on the planar case has been classified by Liu-Lu \cite{liu2022number} for $p=0,q\in(0,4]$, and been classified by Li-Wan \cite{li2022classification} for a more wider range of indices $(p,q)$. Concerning the uniqueness problem for $f$ closing to one, Chen-Feng-Liu \cite{chen2022uniqueness} obtain uniqueness of \eqref{e1.2} on ${\mathbb{R}}^3$  and $p=0, q=n$. A similar result was extended to all dimensions by B\"{o}r\"{o}czky-Saroglou \cite{boroczky2023uniqueness} recently. For the general $f$, Hu-Ivaki \cite{hu2024uniqueness} proved a uniqueness result of \eqref{e1.2} for $p\in(-n-1,-1), q\in[n,n+1)$, using the spectral inequality for convex bodies.

For its importance and difficulty, this is the main purpose for us to discuss the problem. We will show the following multiple solutions result to \eqref{e1.3}.

\begin{theo}\label{t1.1}
 Considering \eqref{e1.3} on dimension $n\geq2$ and assuming $p\leq-n+1$, we have the following result of multiple solutions which are different from each other up to rotations:

 (1) On the planar case $n=2$, if $q-p>\tau^2$ hold for some integer $\tau\geq3$, then \eqref{e1.3} admits at least $\tau-1$ solutions.

(2) When $n=3$, if $q-p>12$, there exist at least two solutions. If $q-p>20$, there exist at least four solutions. If $q-p>42$, there exist at least five solutions.

(3) When $n\geq4$, if $q-p>3(n+1)$, there exists at least two solutions. If $q-p>4(n+2)$, there exist four solutions.
 \end{theo}

\begin{rema} When $n=3, k=20$ or $n=4, k=24, 120, 600$, there exists at least one more $k-$symmetric solution which is not constant ${\mathbb{S}}^{n-1}$, as long as $q-p$ is large enough. Actually, this solution exists for $q-p>\lambda_1(n,k)$, where $\lambda_1(n,k)$ denotes the best constant of Poincar\'{e} inequality on the family of $k-$symmetric functions with zero mean.
\end{rema}

 Our multiple solutions result Theorem \ref{t1.1} improves largely a similar result obtained recently by Chen-Chen-Li (\cite{chen2021variations}, Theorem 1.1). The main ingredient of our proof was based on the study of a variational scheme on the family of functions with large symmetricity. Actually, we have shown the following existence result for symmetric solutions of \eqref{e1.3}.

\begin{theo}\label{t1.2}
  Consider \eqref{e1.3} on dimension $n\geq2$ and assuming $p\leq-n+1$. Supposing that there exists a proper subgroup $\Gamma$ of $SO(n)$ such that $\cup_{\phi\in\Gamma}\phi(e_1), e_1\equiv(1,0,\cdots,0)\in{\mathbb{R}}^n$ spans the whole space ${\mathbb{R}}^n$, then there exists a $\Gamma-$symmetric solution $h$ to \eqref{e1.3} on ${\mathbb{S}}^{n-1}$.
\end{theo}

The second part of this paper is devoted to the discussion of the solvability problem of \eqref{e1.2}. Existence result of weak solution for rather weak assumption on $f$ has been obtained by Zhao \cite{zhao2017dual} for $p=0$ and $q<0$, have been obtained by Huang-Zhao \cite{huang2018lp} for complementary range with respect to $p<0<q$, and have been obtained by Chen-Li \cite{chen2021lp} for negative $p$. The existence result of classical solution for positive H\"{o}lder $f$ can also be found in \cite{huang2018lp} for $p>q$, using the method of continuity. The readers may also refer to the papers \cite{huang2016geometric,zhao2017dual,kim2022diameter,saroglou2021nonexistence,boroczky2022logarithmic,boroczky2023anisotropic,zhao2018existence,guang2023flow,boroczky2018subspace,henk2018necessary,boroczky2012log,lutwak1993brunn,haberl2010even,huang2018lp,zhu2015lp,chen2017lp} for the related topic. In this paper, we will prove the following nonexistence result.

\begin{theo}\label{t1.3}
  Supposing that $n\geq2$ and $p\leq-q, q\geq n$, then there exists a smooth function $f$ which is positive outside equator or outside two polar, such that \eqref{e1.2} has no solution.
\end{theo}

It is also remarkable that in a recent paper by Guang-Li-Wang \cite{guang2022l_p}, solvability of \eqref{e1.2} for positive $f$ and $q=n, p<-n$ has been shown using the method of curvature flow. Therefore, our construction assumption for $f$ vanishing somewhere may not be improved in generally. The proof of Theorem \ref{t1.3} was based on the following Chou-Wang type identity which was derived in \cite{chou2006lp} for $q=n$ and $p=-n$.

\begin{theo}\label{t1.4}
  For any $C^2-$solution $h$ of \eqref{e1.2} with $pq\not=0$ and any tangent projective vector field $\xi$ of the form
    $$
     \xi(x)\equiv Bx-(x^TBx)x, \ \ x\in{\mathbb{S}}^{n-1}
    $$
 for any given $n\times n$ matrix $B$, there holds
  \begin{equation}\label{e1.4}
    \int_{{\mathbb{S}}^{n-1}}\Big(\nabla_\xi f+\frac{p+q}{q}fdiv\xi\Big)h^p=-\frac{p(q-n)}{q}\int_{{\mathbb{S}}^{n-1}}fh^p|Z|^{-2}(h\nabla_\xi h+\nabla\xi(\nabla h,\nabla h)).
  \end{equation}
\end{theo}

This paper was initiated by the first author in the year 2021, who considered the special case when $q=n$ [arXiv:2104.07426,v1] for $L_p$ Minkowski case $q=n$. The general case when $q\neq n$ was recently completed by all the three authors. The content of the paper was organized as follows.  We define the $k-$symmetricity on Section 2. By introducing a new variational scheme on $k-$symmetricity, we prove the existence of positive smooth solution of \eqref{e1.3} which maximizes the variational functional on Section 3. Next, we calculate the second variation for the constant function on Section 4, and then calculate the best constant of Poincar\'{e} inequality upon the family of $k-$symmetric functions on Section 5. As an application, the proof of Theorem \ref{t1.1} was completed as a consequence of Theorems \ref{t4.1} and \ref{t5.1}. After proving a new Chou-Wang type identity \eqref{e1.4} on Section 7, we present our Main Nonexistence Theorem \ref{t1.3} on the final section.

\vspace{10pt}

\section{$k-$symmetricity on a discrete sub-group of $SO(n)$}

In order discussing the uniqueness problem of \eqref{e1.3}, we need to introduce a variational scheme on the range $p\leq-n+1$ and $q> p+\lambda_1(n,k)$. Since the lacking of Blaschke-Santalo's inequality, one needs to restrict the family of functions into a $k-$symmetricity class defined as following.

\begin{defi}\label{d2.1}
 Taking a $k-$regular polytope $\vartheta^{(k)}$ with vertices $\vartheta=(\vartheta_1,\cdots,\vartheta_{k})\in{\mathbb{S}}^{n-1}$ spreading evenly on ${\mathbb{S}}^{n-1}$, we define a subgroup ${\mathcal{S}}_{k}(\vartheta)$ of $SO(n)$ by mappings $\phi\in SO(n)$ which map $\{\vartheta_j\}_{j=1}^{k}$ onto $\{\vartheta_j\}_{j=1}^{k}$. A function $h\in C^\infty({\mathbb{S}}^{n-1})$ (or a function $\xi\in C^\infty({\mathbb{R}}^{n})$ is called to be $k-$symmetric and denoted by $h\in\widetilde{{\mathcal{S}}}_{k}(\vartheta)$ (or $\xi\in\widetilde{{\mathcal{S}}}_{k}(\vartheta)$) if and only if it is invariant under all mappings $\phi\in{\mathcal{S}}_{k}(\vartheta)$.
\end{defi}

The concept of symmetry group has been studied extensively in different topics of differential geometry and differential equations. The readers may refer to \cite{olver1988applications,lunardi2012analytic,andrews2000motion} for some related issues. Concerning the numbers of $k-$regular polytopes with different shapes on ${\mathbb{R}}^n$, we have the following elementary proposition.

\begin{prop}\label{p2.1} We have the following alternatives concerning the number of $k-$regular polytopes with different shapes on Euclidean space ${\mathbb{R}}^{n}$:

(1) When $n=2$, there exist $k-$regular polytopes for each $k\in{\mathbb{N}}, k\geq3$.

(2) When $n=3$, there exist only $k-$regular polytopes for $k=4, 6, 8, 12, 20$.

(3) When $n=4$, there exist only $k-$regular polytopes for $k=5, 8, 16, 24, 120, 600$.

(4) When $n\geq5$, there exist only $k-$regular polytopes for $k=n+1, 2n, 2^n$.
\end{prop}

\noindent This proposition is well known. The explicit construction of $2n$-regular polytope follows from
    $$
    \begin{cases}
     \vartheta^{(2n)}=(\vartheta_1,\cdots,\vartheta_{2n}),\\
      \vartheta_{\pm j}=\pm e_j\in{\mathbb{S}}^{n-1}, \ \ j=1,2,\cdots, n
    \end{cases}
    $$
and the explicit construction of $2^n$-regular polytope follows from
    $$
    \begin{cases}
     \vartheta^{(2^n)}=(\vartheta_1,\cdots,\vartheta_{2n}),\\
      \vartheta_{j}=\frac{(\pm 1,\pm 1,\cdots,\pm 1)}{\sqrt{n}}\in{\mathbb{S}^{n-1}}, \ \ j=1,2,\cdots, 2^n,
    \end{cases}
    $$
where $e_j$ stands for the $j-$th axial unit vector whose $j-$th entry is one and the others are zero. We present here only a simple proof for the existence of $(n+1)$-regular polytopes for each dimension $n\geq2$. When $n=2$, the conclusion is clearly true. Actually, trisection points $\vartheta_1, \vartheta_2, \vartheta_3$ of ${\mathbb{S}}^1$ are desired. Assume the conclusion is true for $n=k$, we show that it is also true for $n=k+1$. Actually, noting that
   $$
    {\mathbb{S}}^{k}\cap\Big\{y\in{\mathbb{R}}^{k+1}\big|\ y_{k+1}=0\Big\}={\mathbb{S}}^{k-1},
   $$
one can choose $k+1$ points
  $$
    \overline{\vartheta}_1, \overline{\vartheta}_2, \cdots \overline{\vartheta}_{k+1}\in{\mathbb{S}}^{k}\cap\Big\{y\in{\mathbb{R}}^{k+1}\big|\ y_{k+1}=0\Big\}, \ \ \overline{\vartheta}_i=(\overline{y}'_i, 0)
  $$
which are spreading evenly on ${\mathbb{S}}^{k-1}$ by induction hypothesis. We denote $l_{k+1}$ the length of the regular polytope determined by $\vartheta_1, \vartheta_2, \cdots, \vartheta_{k+1}$. Then, $k+2$ points $\vartheta_1, \vartheta_2, \cdots, \vartheta_{k+2}$ spreading evenly on ${\mathbb{S}}^{k}$ can be constructed by setting
   $$
    \vartheta_{k+2}=(0,1), \ \ \vartheta_i=\Big(\sqrt{1-z^2}\overline{y}'_i,-z\Big), \ \ \forall i=1,2,\cdots, k+1,
   $$
where $z\equiv\frac{l_{k+1}^2-2}{l_{k+1}^2}$. Actually, direct calculation shows that $\vartheta_1,\cdots, \vartheta_{k+2}$ determine a regular polytope with length
       $$
        l_{k+2}=\frac{2\sqrt{l_{k+1}^2}-1}{l_{k+1}}.
       $$
So, the conclusion follows from mathematical induction. $\Box$

In order clarifying the concept of $k-$symmetricity which will be used to introduce our new variational scheme, we give a proof to the following fundamental proposition.

\begin{prop}\label{p2.2}
 For any $k-$regular polytope $\vartheta^{(k)}$ with vertices $\vartheta=(\vartheta_1,\cdots,\vartheta_k)\in{\mathbb{S}}^{n-1}$ presented in Proposition \ref{p2.1} for $n=2, k\geq3$ or $n\geq3$ and $k=n+1, 2n, 2^n$, the subgroup ${\mathcal{S}}_k(\vartheta)$ defined in Definition \ref{d2.1} is a nonempty finite discrete proper subgroup of $SO(n)$.
\end{prop}

\noindent\textbf{Proof.} On the planar case $n=2$, the conclusion of this proposition is clearly true. Since ${\mathcal{S}}_k(\vartheta)$ is composed exactly of rotations $\phi_j, j=0,1,\cdots, k-1$ with respect to rotation angles $2j\pi/k, j=0,1,\cdots, k-1$. When $n\geq3$ and $k=n+1$, nonempty of ${\mathcal{S}}_k(\vartheta)$ follows from the following construction. At firstly, since each three vertices $\vartheta_{i_1}, \vartheta_{i_2}, \vartheta_{i_3}$ of $n+1$-regular polytope forms a regular triangle. So, it is clearly that there exists a rotation $\phi\in SO(n)$ such
  $$
   \phi(\vartheta_{i_1})=\phi(\vartheta_{i_2}), \ \ \phi(\vartheta_{j})=\phi(\vartheta_j), \ \ \forall j\not=i_1, i_2.
  $$
That is to say, ${\mathcal{S}}_k(\vartheta)$ is not empty for $k=n+1$. When $n\geq3$ and $k=2n$, without loss of generality, one may assume that the vertices of $\vartheta^{(2n)}$ are given by
  $$
   \vartheta_{\pm j}=\pm e_j\in{\mathbb{S}}^{n-1}, \ \ \forall j=1,2,\cdots, n.
  $$
Hence, for each $a\not=b$, one can use the rotation $\phi\in SO(n)$ on the axis plane of $(e_a, e_b)$ which maps one of $\pm e_a$ to another one of $\pm e_b$, and reserves the other vertices $\vartheta_{\pm j}, j\not=a, b$. A similar existence result of $\phi\in SO(n)$ maps $e_i$ to $-e_i$ is also clearly true. In the final case $n\geq3$ and $k=2^n$, one may assume that the vertices of $\vartheta^{(2^n)}$ was formed by
  $$
   \vartheta_{j}=\frac{(\pm 1, \pm 1,\cdots, \pm 1)}{\sqrt{n}}\in{\mathbb{S}}^{n-1}, \ \ \forall j=1,2,\cdots, 2^n.
  $$
Note that any two vertices $\vartheta_{a}, \vartheta_{b}$ can be connected by
  $$
   \vartheta_{a}=\vartheta_{j_1}-\vartheta_{j_2}-\cdots-\vartheta_{j_r}=\vartheta_{b}, \ \ r\geq3
  $$
such that $\vartheta_{j_i}$ and $\vartheta_{j_{i+1}}$ lie on a same face of $\vartheta^{(2^n)}$ for each $i$. Since for each $i$, there exists $\phi\in SO(n)$ maps $\vartheta_{j_i}$ to $\vartheta_{j_{i+1}}$, a composition of these maps gives the desired mapping $\phi\in SO(n)$ maps $\vartheta_a$ to $\vartheta_b$. The nonempty of ${\mathcal{S}}_k(\vartheta^{(k)})$ was shown.

Now, we remain to show that ${\mathcal{S}}_k(\vartheta)$ is a finite subgroup of $SO(n)$. Actually, note first that $\phi(\vartheta^{(k)})=\phi'(\vartheta^{(k)}), \phi,\phi'\in SO(n)$ implies that $\phi=\phi'$, since $\vartheta^{(k)}$ spans a whole space ${\mathbb{R}}^n$. Therefore, ${\mathcal{S}}_k(\vartheta^{(k)})$ is isomorphism to a permutation group of $\vartheta^{(k)}$, and thus is finite. The proof of the proposition was done. $\Box$

\begin{rema}
  For the other cases
  $$
   n=3, k=12, 20, \ \mbox{ or } \ n=4, k=24, 120, 600,
  $$
the conclusion of Proposition \ref{p2.2} is also true. Actually, one can prove by Mathematical Induction. Assuming the nonempty of ${\mathcal{S}}_{k'}(\vartheta^{(k')})$ holds for any $k'$-regular polytope with $k'<k$ on any dimensions $n'<n$ . Then, given any two vertices $\vartheta_{a}, \vartheta_{b}, a\not=i_2b$, if they lie on a same face $\vartheta^{k'}$ of $\vartheta^{(k)}$, the conclusion is clearly true by induction hypothesis, since for $k$-regular polytope, rotation on a face makes $\vartheta^{(k)}$ invariant as long as it makes this face invariant.  If these two vertices does not lie on a same face of $\vartheta^{(k)}$, one can connect them by
  $$
   \vartheta_{a}=\vartheta_{j_1}-\vartheta_{j_2}-\cdots-\vartheta_{j_r}=\vartheta_{b}, \ \ r\geq3
  $$
such that $\vartheta_{j_i}$ and $\vartheta_{j_{i+1}}$ lie on a same face of $\vartheta^{(k)}$ for each $i$. Then the existence of $\phi\in{\mathcal{S}}_k(\vartheta^{(k)})$ follows from composition.
\end{rema}

\vspace{10pt}

\section{New variational scheme on the $k-$symmetricity}

From now on, we assume that there exists $k-$regular polytope with vertices $\vartheta=(\vartheta_1,\cdots,\vartheta_k)\in{\mathbb{S}}^{n-1}$. Letting ${\mathcal{C}}$ be the family of support function $u$ of convex body $\Omega\in{\mathcal{K}}_0$, we set
   $$
    {\mathcal{C}}_k\equiv\Big\{h\in\widetilde{{\mathcal{S}}}_{k}(\vartheta)\cap{\mathcal{C}}\Big|\ \int_{{\mathbb{S}}^{n-1}}h^pdx=\alpha_n\equiv|{\mathbb{S}}^{n-1}|\Big\}
   $$
to be the symmetric class and consider a variational scheme
   \begin{equation}\label{e3.1}
     \sup_{h\in{\mathcal{C}}_k}F(h), \ \ F(h)\equiv \frac{1}{q}\ln\int_{{\mathbb{S}}^{n-1}}r^q(y)dy
   \end{equation}
on ${\mathcal{C}}_k$. We have the following {\it a-priori} upper bound for functions $h\in{\mathcal{C}}_k$.

\begin{lemm}\label{l3.1}
  Letting $n\geq2$ and $p<0$, there exists a positive constant $C_{n,p}$ depending only on $n, p$, such that
     \begin{equation}\label{e3.2}
        h(x)\leq C_{n,p}, \ \ \forall x\in{\mathbb{S}}^{n-1}
     \end{equation}
  holds for each $h\in{\mathcal{C}}_k$.
\end{lemm}

\noindent\textbf{Proof.} Noting that if the maximum of a $k-$symmetric function $h$ becomes large, the minimum of $h$ must tend to infinity as well by convexity. Then, it is contradicting with the fact that $\alpha_n=\int_{{\mathbb{S}}^{n-1}}h^p$. The proof of the lemma was done. $\Box$\\

The next lemma estimates the functions in ${{\mathcal{C}}_k}$ from below.

\begin{lemm}\label{l3.2}
 Letting $n\geq2$ and $p\leq-n+1$, there exists a positive constant $C_{n,p}$ depending only on $n, p$, such that
     \begin{equation}\label{e3.3}
        h(x)\geq C^{-1}_{n,p}, \ \ \forall x\in{\mathbb{S}}^{n-1}
     \end{equation}
  holds for any function $h\in{{\mathcal{C}}_k}$.
\end{lemm}

\noindent\textbf{Proof.} Without loss of generality, one may assume that the minimum of $h$ is attained at south polar of ${\mathbb{S}}^{n-1}$. Projecting the south hemisphere ${\mathbb{S}}^{n-1}_-\equiv\big\{x=(x',x_{n})\in{\mathbb{S}}^{n-1}\big|\ x_{n}<0\big\}$ onto the hyperplane $x_{n}=-1$ by the mapping
   $$
    T(x)=z, \ \ x\equiv\Big(\frac{z}{\sqrt{1+|z|^2}},\frac{-1}{\sqrt{1+|z|^2}}\Big), \ \ z\in{\mathbb{R}}^{n-1}
   $$
and setting $v(y)=\sqrt{1+|y|^2}h(x)$, one gets that $D^2_{ij}v(y)=\sqrt{1+|y|^2}(\nabla^2_{ij}h+h\delta_{ij})$ is positive definite. Thanks to the boundedness and convexity of $v$ on $B_2$, the Lipschitz bound of $v$ on $B_1$ can not be greater than $osc_{B_2}(v)$. So, letting $m\equiv h(e)=\min_{x\in{\mathbb{S}}^{n-1}}h(x)$ for $e=(0,-1)\in{\mathbb{S}}^{n-1}_-$, it is inferred that
  \begin{equation}\label{e3.4}
       |v(y)-m|\leq |y|osc_{B_2}(v)\leq C_1|y|, \ \ \forall y\in B_1
   \end{equation}
holds for some positive constant $C_1$ depending only on $n$ and $p$. Combining with
   \begin{equation}\label{e3.5}
      \int_{B_1}v^pdy\leq C\int_{{\mathbb{S}}^{n-1}}h^p\leq C_2
   \end{equation}
for some positive constant $C_2$ depending only on $n, p$, we conclude that
   \begin{eqnarray*}
     C_2&\geq&\int_{B_1}v^pdy\geq\int_{B_1}(m+C_1|y|)^p\\
      &\geq&\begin{cases}
        C_3^{-1}m^{p+n-1}-C_4, & p<-n+1,\\
        -C_3^{-1}\log m-C_4, & p=-n+1
      \end{cases}
   \end{eqnarray*}
for some positive constants $C_3, C_4$ depending only on $n$ and $p$. So, the inequality \eqref{e3.3} follows and the proof of the lemma was done. $\Box$\\

By Lemmas \ref{l3.1} \ref{l3.2}, the functional $F$ is bounded from both below and above. Taking a maximizing sequence of convex bodies $\Omega_j$ with respect to $u_j\in{{\mathcal{C}}_k}, j\in{\mathbb{N}}$ of $F$, it follows from Blaschke's selection theorem that for a subsequence, $\{\Omega_j\}$ converges to a limiting convex body $\Omega_\infty$ with positive Lipschitz support function $h_\infty\in{\mathcal{C}}_k$. After a slightly modifying, by an argument by Chen-Chen-Li in \cite{chen2021variations}, the limiting function $h_\infty$ satisfies the Euler-Lagrange equation \eqref{e1.3} up to a constant. Then, $C^0$-regularity implies $C^1-$regularity by
  $$
   |Z_\infty|^2=h_\infty^2+|\nabla h_\infty|^2\leq\max_{x\in{\mathbb{S}}^{n-1}}h_\infty^2, \ \ \forall y\in{\mathbb{S}}^{n-1}.
  $$
Henceforth, an application of the $C^{1,\alpha}$ estimation in \cite{caffarelli1990localization} and then $C^{2,\alpha}$ estimation in \cite{caffarelli1990interior} by Caffarelli yields the following existence result.

\begin{theo}\label{t3.1}
  For $n\geq2, q\not=0$ and $p\leq-n+1$, there exists a smooth positive maximizer $h_\infty$ of variational problem \eqref{e3.1} satisfying the Euler-Lagrange equation \eqref{e1.3} up to a constant. Moreover, the solution $h_\infty$ can not be a constant unless it is identical to one.
\end{theo}

\vspace{10pt}

\section{Instability of constant solution for $q>p+\lambda_1(n,k)$}

In this section, we discuss the stability of constant function $h\equiv1$. A critical point $h$ of $F(\cdot)$ is called to be stable, if the second variation of the functional at $h$ is non-positive. Otherwise, we will call it to be unstable. We have the following theorem.

\begin{theo}\label{t4.1}
  Supposing that $n\geq2$ and there exists a $k-$regular polytope $\vartheta^{(k)}$, if
     \begin{equation}\label{e4.1}
       q>p+{\lambda_1(n,k)}
     \end{equation}
holds for the first eigenvalue of $-\triangle_{{\mathbb{S}}^{n-1}}$ on the family of $k-$symmetricity functions with zero mean,
 the constant solution  $h\equiv 1$ is a critical point of the functional $F$ which is unstable. As a result, $h\equiv1$ is not a maximizer of $F$ on the family ${{\mathcal{C}}_k}$.
\end{theo}

\begin{rema}
 The best constant ${{\lambda_1(n,k)}}$ of $-\triangle_{{\mathbb{S}}^{n-1}}$ depends only only on $k$ but not on the position of the vertices. Actually, for each pair of $k-$regular polytopes $\vartheta^{(k)}$ and $\widetilde{\vartheta^{(k)}}$, there exists $\phi\in SO(n)$ such that
    $$
     \widetilde{\vartheta^{(k)}}=\phi(\vartheta^{(k)}).
    $$
Since the Rayleigh quotient of Poincar\'{e} inequality in invariant under rotation of ${\mathbb{R}}^{n}$, it is inferred that the best constants $\lambda_1(n,k)$ of $\vartheta^{(k)}$ and $\widetilde{\vartheta^{(k)}}$ are the same.
\end{rema}

To show the conclusion of Theorem \ref{t4.1}, let us calculate the second variation of $F$ in the following lemma. The second variational formula \eqref{e4.3} under below was firstly calculated by Chen-Chen-Li (\cite{chen2021variations}, Proposition 2.1). We present here a simplified proof for the convenience of the readers.
  \begin{lemm}\label{l4.1}
    Under the normalized variation
   \begin{equation}\label{e4.2}
    \varphi_\varepsilon\equiv\varphi(\varepsilon,\theta)\equiv\frac{1+\varepsilon\xi}{\big(\fint_{{\mathbb{S}}^{n-1}}(1+\varepsilon\xi)^p\big)^{1/p}}, \ \ \xi\in\widetilde{{\mathcal{S}}_k}(\vartheta)
   \end{equation}
  near $h=1$, we have the second variational formula
   \begin{equation}\label{e4.3}
     \alpha_n\frac{d^2}{d\varepsilon^2}\Big|_{\varepsilon=0}F(\varphi_\varepsilon)=(q-p)\int_{{\mathbb{S}}^{n-1}}|\xi-\overline{\xi}|^2-\int_{{\mathbb{S}}^{n-1}}|\nabla\xi|^2,
   \end{equation}
  where $\overline{\xi}\equiv\alpha_n^{-1}\int_{{\mathbb{S}}^{n-1}}\xi$ and
   $$
    F(\varphi_\varepsilon)=\frac{1}{q}\ln\int_{{\mathbb{S}}^{n-1}}r^{q-n-1}_\varepsilon\varphi_\varepsilon\det(\nabla^2\varphi_\varepsilon+\varphi_\varepsilon I),\ \ r_\varepsilon^2=\varphi_\varepsilon^2+|\nabla\varphi_\varepsilon|^2.
   $$
  \end{lemm}

\noindent\textbf{Proof.} Direct computation shows that
   \begin{eqnarray}\nonumber\label{e4.4}
     \varphi'_\varepsilon&\equiv&\frac{d}{d\varepsilon}\Big|_{\varepsilon=0}\varphi(\varepsilon,\theta)=\xi-\alpha_n^{-1}\int_{{\mathbb{S}}^{n-1}}\xi\\
    \ \ \ \ \ \ \ \  \varphi''_\varepsilon&\equiv&\frac{d^2}{d\varepsilon^2}\Big|_{\varepsilon=0}\varphi(\varepsilon,\theta)=-2\alpha_n^{-1}\xi\int_{{\mathbb{S}}^{n-1}}\xi\\ \nonumber
     &&-(p-1)\alpha_n^{-1}\int_{{\mathbb{S}}^{n-1}}\xi^2+(p+1)\alpha_n^{-2}\Big(\int_{{\mathbb{S}}^{n-1}}\xi\Big)^2.
   \end{eqnarray}
Thus, we derive
   \begin{equation}\label{e4.5}
     q\alpha_n\frac{d^2}{d\varepsilon^2}\Big|_{\varepsilon=0}F(\varphi_\varepsilon)=\frac{d^2}{d\varepsilon}\Big|_{\varepsilon=0}G_\varepsilon-\alpha_n^{-1}\Big(\frac{d}{d\varepsilon}\Big|_{\varepsilon=0}G_\varepsilon\Big)^2
   \end{equation}
for
  $$
   G_\varepsilon\equiv\int_{{\mathbb{S}}^{n-1}}r^{q-n}_\varepsilon\varphi_\varepsilon\det(\nabla^2\varphi_\varepsilon+\varphi_\varepsilon I).
  $$
Note that
   \begin{equation}\label{e4.6}
    \frac{d}{d\varepsilon}\Big|_{\varepsilon=0}G_\varepsilon=\int_{{\mathbb{S}}^{n-1}}\varphi'_\varepsilon+(q-n)\int_{{\mathbb{S}}^{n-1}}r'_\varepsilon+\int_{{\mathbb{S}}^{n-1}}U^{ij}(\nabla^2_{ij}\varphi'_\varepsilon+\varphi'_\varepsilon\delta_{ij})
   \end{equation}
and
   \begin{eqnarray}\nonumber\label{e4.7}
    \frac{d^2}{d\varepsilon^2}\Big|_{\varepsilon=0}G_\varepsilon&=&\int_{{\mathbb{S}}^{n-1}}\varphi''_\varepsilon+(q-n)\int_{{\mathbb{S}}^{n-1}}\Big[r''_\varepsilon+(q-n-1)r'^2_\varepsilon\Big]+\int_{{\mathbb{S}}^{n-1}}U^{ij}(\nabla^2_{ij}\varphi''_\varepsilon+\varphi''_\varepsilon\delta_{ij})\\
    &&+\int_{{\mathbb{S}}^{n-1}}U^{ij,rs}(\nabla^2_{ij}\varphi'_\varepsilon+\varphi'_\varepsilon\delta_{ij})(\nabla^2_{rs}\varphi'_\varepsilon+\varphi'_\varepsilon\delta_{rs})+2\int_{{\mathbb{S}}^{n-1}}\varphi'_\varepsilon U^{ij}(\nabla^2_{ij}\varphi'_\varepsilon+\varphi'_\varepsilon\delta_{ij})\\ \nonumber
    &&+2(q-n)\int_{{\mathbb{S}}^{n-1}}r'_\varepsilon\varphi'_\varepsilon+2(q-n)\int_{{\mathbb{S}}^{n-1}}r'_\varepsilon U^{ij}(\nabla^2_{ij}\varphi'_\varepsilon+\varphi'_\varepsilon\delta_{ij})
   \end{eqnarray}
hold for
  \begin{equation}\label{e4.8}
    U^{ij}=\delta_{ij}, \ \ U^{ij,rs}=\delta_{ij}\delta_{rs}-\delta_{ir}\delta_{js}
  \end{equation}
and
   \begin{eqnarray}\nonumber\label{e4.9}
    r'_\varepsilon&\equiv&\frac{d}{d\varepsilon}\Big|_{\varepsilon=0}r_\varepsilon=\varphi'_\varepsilon\\
    r''_\varepsilon&\equiv&\frac{d^2}{d\varepsilon^2}\Big|_{\varepsilon=0}r_\varepsilon=\varphi''_\varepsilon+|\nabla\varphi'_\varepsilon|^2.
   \end{eqnarray}
Substituting \eqref{e4.8}-\eqref{e4.9} into \eqref{e4.6}-\eqref{e4.7}, it yields that
  \begin{equation}\label{e4.10}
    \frac{d}{d\varepsilon}\Big|_{\varepsilon=0}G_\varepsilon=\int_{{\mathbb{S}}^{n-1}}\triangle\xi=0
  \end{equation}
and
  \begin{eqnarray}\nonumber\label{e4.11}
   \frac{d^2}{d\varepsilon^2}\Big|_{\varepsilon=0}G_\varepsilon&=&q\int_{{\mathbb{S}}^{n-1}}\varphi''_\varepsilon-(q+n-2)\int_{{\mathbb{S}}^{n-1}}|\nabla\varphi'_\varepsilon|^2+q(q-1)\int_{{\mathbb{S}}^{n-1}}\varphi'^2_\varepsilon\\
   \nonumber   &&-\int_{{\mathbb{S}}^{n-1}}|\nabla^2\varphi'_\varepsilon|^2+\int_{{\mathbb{S}}^{n-1}}|\triangle\varphi'_\varepsilon|^2\\ &=&q(q-p)\int_{{\mathbb{S}}^{n-1}}|\xi-\overline{\xi}|^2-(q+n-2)\int_{{\mathbb{S}}^{n-1}}|\nabla\xi|^2-\int_{{\mathbb{S}}^{n-1}}|\nabla^2\xi|^2+\int_{{\mathbb{S}}^{n-1}}|\triangle\xi|^2\\ \nonumber
   &=&q(q-p)\int_{{\mathbb{S}}^{n-1}}|\xi-\overline{\xi}|^2-q\int_{{\mathbb{S}}^{n-1}}|\nabla\xi|^2,
  \end{eqnarray}
where $\overline{\xi}\equiv\alpha_n^{-1}\int_{{\mathbb{S}}^{n-1}}\xi$. The proof of Lemma \ref{l4.1} was done. $\Box$

Then, Theorem \ref{t4.1} was a direct consequence of the following proposition.

\begin{prop}\label{p4.1}
  Supposing that the Poincar\'{e} inequality
   \begin{equation}\label{e4.12}
     \int_{{\mathbb{S}}^{n-1}}|\nabla\xi|^2\geq{\lambda_1(n,k)}\int_{{\mathbb{S}}^{n-1}}\xi^2
   \end{equation}
holds on the family
   \begin{equation}\label{e4.13}
     \xi\in{{\mathcal{G}}_{n,k}}\equiv\Big\{\xi\in H^1({\mathbb{S}}^{n-1})\Big|\ \int_{{\mathbb{S}}^{n-1}}\xi=0,\ \xi\mbox{ is } k-\mbox{symmetric}\Big\}
   \end{equation}
with the best constant ${\lambda_1(n,k)}>0$, we have the instability of the constant solution $h\equiv1$ for
   \begin{equation}\label{e4.14}
     q>p+{\lambda_1(n,k)}.
   \end{equation}
\end{prop}

\begin{rema}
 The best constant of Poincar\'{e} inequality \eqref{e4.12} on the family \eqref{e4.13} is attained by utilizing a usual variational method for lower semi-continuity of weakly convergence.
\end{rema}

\noindent\textbf{Proof of Proposition \ref{p4.1}.} Letting $\xi$ be the eigenfunction of $-\triangle_{{\mathbb{S}}^{n-1}}$ corresponding to the least eigenvalue ${\lambda_1(n,k)}$, one has
   \begin{equation}\label{e4.15}
     \int_{{\mathbb{S}}^{n-1}}|\nabla\xi|^2={\lambda_1(n,k)}\int_{{\mathbb{S}}^{n-1}}\xi^2.
   \end{equation}
Substituting into \eqref{e4.3} yields that
    $$
     \alpha_n\frac{d^2}{d\varepsilon^2}\Big|_{\varepsilon=0}F(\varphi_\varepsilon)=(q-p-{\lambda_1(n,k)})\int_{{\mathbb{S}}^{n-1}}\xi^2>0
    $$
in case $q>p+{\lambda_1(n,k)}$. The instability of constant solution $u\equiv1$ has been shown. $\Box$

\vspace{10pt}

\section{Best constant of Poincar\'{e} inequality on ${{\mathcal{G}}_{n,k}}$}

\begin{theo}\label{t5.1}
   Supposing that there exists $k$-regular polytope $\vartheta^{(k)}$ on the Euclidean space ${\mathbb{R}}^{n}$ for some $n\geq2, k\geq3$, the best constant of Poincar\'{e}'s inequality \eqref{e4.12} on the family ${{\mathcal{G}}_{n,k}}$ is given by
   \begin{equation}\label{e5.1}
     \lambda_1(n,k):\begin{cases}
        =k^2, &  n=2, \forall k\geq3,\\
        =3(n+1), & n\geq3, k=n+1,\\
        =4(n+2), & n\geq3, k=2n \mbox{ or } 2^{n},\\
        \in[12,42], & n=3, k=12,\\
        \in[12,\infty), & n=3, k=20\\
        \in[15,\infty), & n=4, k=24, 120, 600.
     \end{cases}
   \end{equation}
\end{theo}

\begin{rema}
 The upper bound of $\lambda_1(n,k)$ for $n=3, k=12$ or $n\geq3, k=n+1, 2^{n}$ can be found in \cite{andrews2000motion}. Here, we give a proof to the precise formula \eqref{e5.1} for $\lambda_1(n,k)$ as follows.
\end{rema}

\noindent\textbf{Proof of Theorem \ref{t5.1}.} In the planar case $n=2$, \eqref{e4.12}-\eqref{e4.13} change to
    \begin{equation}\label{e5.2}
      \int^{2\pi}_0\xi_\theta^2\geq\lambda_1(2,k)\int^{2\pi}_0\xi^2
    \end{equation}
and
    \begin{equation}\label{e5.3}
      \xi\in{\mathcal{G}}_{2,k}\equiv\Big\{\xi\in H^1({\mathbb{S}}^1)\Big| \ \int_{{\mathbb{S}}^1}\xi=0, \ \xi\mbox{ is a } 2\pi/k-\mbox{periodic function}\Big\}.
    \end{equation}
Therefore, $\xi\in{\mathcal{G}}_{2,k}$ is equivalent to
   $$
    \xi(\theta)=\Sigma_{j=1}^\infty \Big(a_k\cos(jk\theta)+b_k\sin(jk\theta)\Big),
   $$
which implies that ${\lambda_1(n,k)}=k^2$ for $n=2$.

 When $n\geq3$, the determining of the best constant ${\lambda_1(n,k)}$ of Poincar\'{e} inequality \eqref{e4.12}-\eqref{e4.13} is related to the eigenvalues problem of Laplace-Beltrami operator $-\triangle_{{\mathbb{S}}^{n-1}}$ on ${\mathbb{S}}^{n-1}$. As well known that all eigenfunctions $u$ of $-\triangle_{{\mathbb{S}}^{n-1}}$ are given by homogeneous harmonic polynomials $v(y)=|y|^\mu u\Big(\frac{y}{|y|}\Big)$ of degree $\mu$ on ${\mathbb{R}}^n$, whose eigenvalues are given exactly by $\lambda\equiv\mu(n+\mu-2)$ for positive integer $\mu$. In a paper of Kazdan (\cite{kazdan1998solving}, page 12), the multiplicities of the eigenvalues were also calculated explicitly. Note that the mean value of $u$ equals to zero if and only if
   \begin{equation}\label{e5.4}
     \int_{B_1}v(y)dy=0,
   \end{equation}
and $u$ is $k-$symmetric if and only if $v$ is $k-$symmetric. Setting
   $$
    v(y)=\Sigma_{|\alpha|=\mu}a^\alpha y^\alpha, \ \ \alpha\equiv(\alpha_1,\cdots,\alpha_{n})\in{\mathbb{R}}^{n}
   $$
for $y^\alpha\equiv\Pi_{i=1}^{n}y_i^{\alpha_i}$ as usual, and summing the above relations between $u$ and $v$, we obtain that
   \begin{eqnarray}\nonumber\label{e5.5}
     &&\Sigma_{j=1}^{n}\Sigma_{|\beta|=\mu-2}a^{\beta+2e_j}(\beta_j+2)(\beta_j+1)=0,\\
     &&\Sigma_{|\alpha|=\mu}a^\alpha\int_{B_1}y^\alpha dy=0
   \end{eqnarray}
and
   \begin{equation}\label{e5.6}
    \Sigma_{|\alpha|=\mu}a^\alpha(\phi(y))^\alpha=\Sigma_{|\alpha|=\mu}a^\alpha y^\alpha, \ \ \forall y\in{\mathbb{R}}^{n}, \phi\in{\mathcal{S}}_{k}(\vartheta),
   \end{equation}
where $e_j$ is a vector whose $j-$th coordinate is one and the others are zero.

\begin{prop}\label{p5.1}
  Under the assumptions of Theorem \ref{t5.1}, the set
    \begin{eqnarray}\nonumber\label{e5.7}
      &{{\mathcal{H}}_{n,k}}\equiv\Big\{v\in\widetilde{{\mathcal{S}}}_{k}(\vartheta)\setminus\{0\}\big|\ \int_{B_1}v(y)dy=0, \ v\mbox{ is a}&\\
      &\mbox{homogeneous harmonic polynomial on } {\mathbb{R}}^{n}\Big\}&
    \end{eqnarray}
  is nonempty. Moreover, the least degree of $v\in{{\mathcal{H}}_{n,k}}$ is given by
     \begin{equation}\label{e5.8}
       {\mu_1(n,k)}:\begin{cases}
          =k, &  \mbox{ if } n=2, k\geq3\\
          =3, & \mbox{ if } n\geq3, k=n+1,\\
         \in[3,\infty), & \mbox{ if } n\geq3, k>n+1.
       \end{cases}
     \end{equation}
  As a corollary, $\lambda_1(n,n+1)=3(n+1), \forall n\geq2$.
\end{prop}

\noindent\textbf{Proof.} At first, one knows that the best constant ${{\lambda_1(n,k)}}$ of Poincar\'{e} inequality \eqref{e4.12}-\eqref{e4.13} is attained at some function $u\in{\mathcal{G}}_n$. After expanding $u$ to be a homogeneous harmonic function $v$ of degree ${\mu_1(n,k)}$ determined by
   \begin{equation}\label{e5.9}
     {\mu_1(n,k)}(n+{\mu_1(n,k)}-2)-{{\lambda_1(n,k)}}=0,
   \end{equation}
Liouville property shows that $v$ must be a polynomial, also be $k-$symmetric as well. So, the set ${{\mathcal{H}}_{n,k}}$ is nonempty, and it remains to calculate ${\mu_1(n,k)}$.

In the case $n=2$, it is not hard to see that $\mu_1(2,k)=k, \forall k\geq3$ and all homogenous harmonic polynomials $v\in{\mathcal{H}}_{2,k}$ of degree $k$ are given by
   \begin{eqnarray*}
     v(x,y)&=&A(y^k-C_k^2x^2y^{k-2}+C_k^4x^4y^{k-4}+\cdots)\\
      &&+B(C_k^1xy^{k-1}-C_k^3x^3y^{k-3}+C_k^5x^5y^{k-5}+\cdots)
   \end{eqnarray*}
for each given constants $A$ and $B$. We claim that $v$ is $k-$symmetric for each $A, B$. Actually, noting that
  \begin{eqnarray*}
   v_a(x,y)&\equiv& y^k-C_k^2x^2y^{k-2}+C_k^4x^4y^{k-4}+\cdots\\
    &=&(y+xi)^k+(y-xi)^k\equiv z^k+\overline{z}^k,\\
   v_b(x,y)&\equiv& i\Big[C_k^1xy^{k-1}-C_k^3x^3y^{k-3}+C_k^5x^5y^{k-5}+\cdots\\
   &=&(xi+y)^k-(xi-y)^k\equiv z^k-(-\overline{z})^k
  \end{eqnarray*}
for complex number $z\equiv y+xi$, one can conclude the $k-$symmetricity of $v_a, v_b$.

For dimension $n\geq3$, the proof of Proposition \ref{p5.1} was decomposed into the following two lemmas.

\begin{lemm}\label{l5.1}
  For $n\geq3$ and some $(n+1)-$regular polytope $\vartheta^{(n+1)}$ with vertices $\vartheta=(\vartheta_1,\cdots,\vartheta_{n+1})\in{\mathbb{S}}^{n-1}$, there exists at least one homogeneous harmonic polynomial $v$ of degree $3$ which is $(n+1)-$symmetric and satisfies the property
     \begin{equation}\label{e5.10}
       \int_{B_1}v(y)dy=0
     \end{equation}
  of zero mean.
\end{lemm}

\noindent\textbf{Proof.} Taking $k$ linear functions
   \begin{equation}\label{e5.11}
     l_j(y)\equiv \langle \vartheta_j,y\rangle, \ \ \forall j=1,\cdots,n+1, \ y\in{\mathbb{R}}^{n},
   \end{equation}
and setting
   $$
     v_{3}(y)\equiv\Sigma_{i<j<m}l_i(y)l_j(y)l_m(y), \ \ \forall y\in{\mathbb{R}}^{n},
   $$
it is clear that $v_3$ is a $(n+1)-$symmetric polynomial of homogeneous degree $3$. We claim that $v_{3}$ is harmonic and has zero mean in sense of \eqref{e5.10} over the ball $B_1$. Actually, the harmonicity of $v_{3}$ follows from
   $$
    \triangle v_{3}=\frac{n(n-1)}{4}\Sigma_{j=1}^{n+1}l_j(y)=\frac{n(n-1)}{4}\langle\Sigma_{j=1}^{n+1}\vartheta_j,y\rangle=0
   $$
To show the zero mean property, after restricting to ${\mathbb{S}}^{n-1}$, the function
   $$
    u(x)\equiv v_3(x), \ \ \forall x\in{\mathbb{S}}^{n-1}
   $$
satisfies the eigen-equation
   \begin{equation}\label{e5.12}
     -\triangle_{{\mathbb{S}}^{n-1}}u=3(n+1)u, \ \ \forall x\in{\mathbb{S}}^{n-1}.
   \end{equation}
Integrating over the sphere, it yields that
   \begin{eqnarray*}
     \int_{B_1}v_{3}(y)=\frac{1}{n+3}\int_{{\mathbb{S}}^{n-1}}u=0.
   \end{eqnarray*}
So, the claim holds true and the proof of the lemma was done. $\Box$\\

\begin{lemm}\label{l5.2}
  For $n\geq3$ and some $k-$regular polytope $\vartheta^{(k)}$ with vertices $\vartheta=(\vartheta_1,\cdots,\vartheta_k)\in{\mathbb{S}}^{n-1}$, the degree of polynomial $v\in{{\mathcal{H}}_{n,k}}$ can not be less than $3$.
\end{lemm}

\noindent\textbf{Proof.} Supposing for some $a=\Sigma_{j=1}^{k}\varsigma_j\vartheta_j\in{\mathbb{R}}^{n}$,
   $$
    v(y)\equiv\langle a,y\rangle=\Sigma_{j=1}^{k}\varsigma_j\langle \vartheta_j,y\rangle
   $$
is a homogeneous harmonic polynomial of degree one. By $k-$symmetricity of $v$, one has
   $$
    v(y)=C_{n,k}\Sigma_{j=1}^k\langle\Sigma_{i=1}^{k}\vartheta_i,y\rangle=0
   $$
for some constant $C_{n,k}$ depending only on $n$ and $k$. This contradicts with the non-triviality of $v$. Hence, the degree of polynomial $v\in{{\mathcal{H}}_{n,k}}$ can not be one. It remains to show that it can also not be two. Otherwise, let
   $$
    v(y)=\Sigma_{i,j=1}^{n}a_{ij}y_iy_j\in{{\mathcal{H}}_{n,k}}.
   $$
After rotating the coordinates, one may assume that
   \begin{equation}\label{e5.13}
     a_{ij}=\lambda_i\delta_{ij}, \ \ \forall i,j=1,2,\cdots,n.
   \end{equation}
Fixing $i,j=1,\cdots,k, i\not=j$ and assuming $\vartheta_a\equiv(\vartheta_{a1},\cdots,\vartheta_{an})$ for $a=1,2,\cdots, k$, we have by $k-$symmetricity of $v$
   \begin{equation}\label{e5.14}
      \Sigma_{a=1}^{n}\lambda_a\vartheta_{ia}^2=\Sigma_{a=1}^{n}\lambda_a\vartheta_{ja}^2.
   \end{equation}
Introducing an auxiliary variable $d$, we can transform \eqref{e5.14} into the following system
   \begin{equation}\label{e5.15}
      \Sigma_{a=1}^{n}\lambda_a\vartheta_{ja}^2+d=0, \ \ \forall j=1,2,\cdots, k.
   \end{equation}
Since $\{\vartheta_j\}_{j=1}^{k}$ spread evenly on ${\mathbb{S}}^{n-1}$, the rank of the $k\times(n+1)$ coefficient matrix
   $$
    A\equiv\left(
      \begin{array}{ccccc}
        \vartheta_{11}^2 & \vartheta_{12}^2 & \cdots & \vartheta_{1(n+1)n}^2 & 1\\
        \vartheta_{21}^2 & \vartheta_{22}^2 & \cdots & \vartheta_{2n}^2 & 1\\
        \cdots & \cdots & \cdots& \cdots & \cdots \\
        \vartheta_{k1}^2 & \vartheta_{k2}^2 & \cdots & \vartheta_{kn}^2 & 1\\
      \end{array}
    \right)
   $$
is no less than $n$. As a result, the dimension of the set of the solutions equals to one. Combining with the fact that
   $$
     (\lambda_1,\lambda_2,\cdots,\lambda_{n},d)=(1,1,\cdots,1,-1)
   $$
is a non-trivial solution of \eqref{e5.15}, one gets that
   \begin{equation}\label{e5.16}
     \lambda_i=\kappa=-d, \ \ \forall i=1,2,\cdots,n
   \end{equation}
give all solutions for constant $\kappa$. However, it is inferred from harmonicity of $h$ that
   $$
    \triangle_{{\mathbb{R}}^{n}}v(y)=2n\kappa=0.
   $$
So, we derived the desired conclusion $h\equiv0$ and completed the proof of this lemma. $\Box$\\

Now, Proposition \ref{p5.1} is exactly a corollary of Lemma \ref{l5.1} and \ref{l5.2}. To continue the proof of Theorem \ref{t6.1}, we need the following proposition.

\begin{prop}\label{p5.2}
  Supposing that there exists a $k-$polytope $\vartheta^{(k)}$ for dimension $n\geq3$ and some $k=2n$ or $k=2^n$, we have $\lambda_1(n,k)=4(n+2)$.
\end{prop}

The proof of Proposition \ref{p5.2} was divided by the three following lemmas.

\begin{lemm}\label{l5.3}
  Letting $\vartheta^{(2n)}=(\vartheta_1,\cdots,\vartheta_{2n})\in{\mathbb{S}}^{n-1}$ be the regular cross polytope on ${\mathbb{R}}^{n}$ for $n\geq3$, we have the best constant of Poincar\'{e} inequality satisfies
    \begin{equation}\label{e5.17}
      \lambda_1(n,2n)\in\{3(n+1), 4(n+2)\}.
    \end{equation}
\end{lemm}

\noindent\textbf{Proof.} After rotating, one may assume that
  $$
   \vartheta_1=(1,0, \cdots, 0), \ \ \vartheta_2=(-1,0, \cdots, 0)\cdots, \ \ \vartheta_{2n-1}=(0,0, \cdots, 0, 1), \ \ \vartheta_{2n}=(0,0, \cdots, 0, -1).
  $$
Noting that
  $$
    \varphi_3(x,y,z)\equiv x^2y^2+y^2z^2+z^2x^2-\frac{x^4+y^4+z^4}{3}
  $$
is a harmonic function which is symmetric with respect to
$$
 (\pm1,0,0), \ \ (0,\pm1,0), \ \ (0, 0, \pm1),
$$
 if one sets
  $$
   u(x)=\Sigma_{1\leq i<j<k\leq n}\varphi_3(x_i,x_j, x_k),
  $$
we have $u$ is a harmonic function on ${\mathcal{G}}_{n,2n}$ of order four. So, we conclude that $\lambda_1(n,2n)\leq4(n+2)$ as above. $\Box$\\

\begin{lemm}\label{l5.4}
  Letting $\vartheta^{(2^{n})}=(\vartheta_1,\cdots,\vartheta_{2^{n}})\in{\mathbb{S}}^{n-1}$ be a regular hypercube on ${\mathbb{R}}^{n}$ for $n\geq2$, we have the best constant of Poincar\'{e} inequality satisfies
    \begin{equation}\label{e5.18}
      \lambda_1(n,2^{n})\in\{3(n+1), 4(n+2)\}.
    \end{equation}
\end{lemm}

\noindent\textbf{Proof.} After rotating, one may assume that the $j$th entry of the vertex $\vartheta_{i}$ equals to $\pm 1$ for each $i=1,2,\cdots,2^n$ and $j=1,2,\cdots, n$. Noting that $\varphi_2(x,y)\equiv x^3y-xy^3$ is a $4-$symmetric harmonic function on ${\mathbb{R}}^2$, if one sets
 $$
   u(x)=\Sigma_{1\leq i<j\leq n}\varphi_2(x_i,x_j),\ \ x=(x_1,\cdots, x_{n})\in{\mathbb{R}}^{n},
 $$
we have $u$ is a harmonic function on ${\mathcal{G}}_{n,2^{n}}$ of order four. So, we conclude that $\lambda_1(n,2^{n})\leq4(n+2)$ as above. $\Box$\\

\begin{lemm}\label{l5.5}
  Letting  $\vartheta^{(k)}=(\vartheta_1,\cdots,\vartheta_{k})\in{\mathbb{S}}^{n-1}$ be a regular $k-$polytope on ${\mathbb{R}}^{n}$ for $n\geq3$ and $k=2n$ or $k=2^{n}$, we have the best constant of Poincar\'{e} inequality satisfies
    \begin{equation}\label{e5.19}
      \lambda_1(n,k)\not=3(n+1), \ \ k=2n, 2^{n}.
    \end{equation}
\end{lemm}

\noindent\textbf{Proof.} Consider first the case $k=2n$. Note that after rotating, $\vartheta^{(k)}$ can be generated from
  \begin{eqnarray*}
   e&\equiv&(e_1,e_2,\cdots,e_6),\ \ \ \  e_1\equiv(1,0,0), e_2\equiv(-1,0,0),\\
  e_3&\equiv&(0,1,0), e_4\equiv(0,-1,0), e_5\equiv(0,0,1), e_6\equiv(0,0,-1)
  \end{eqnarray*}
and the replacing mapping $\phi_{p,q,r}$ which maps $(x_1,x_2,x_3)$ to $(x_p, x_q, x_r)$ for $p<q<r$. Thus, for each $2n-$symmetric polynomial $u$, there holds
   \begin{equation}\label{e5.20}
    u(x)=\frac{1}{6}\Sigma_{1\leq p<q<r\leq n}\Sigma_{\phi_\alpha\in{\mathcal{S}}_6(e)}\phi_{p,q,r}\circ\phi_\alpha\circ u_3((x_1, x_2, x_3)
   \end{equation}
for $u_3(x_1,x_2,x_3)$ being the sum of terms of $u$ containing only the variables $(x_1,x_2,x_3)$. We claim that for $u$ being polynomial of degree three,
   \begin{equation}\label{e5.21}
     \Sigma_{\phi_\alpha\in{\mathcal{S}_6(e)}}u\big(\phi_\alpha(x_1, x_2, x_3)\big)=0.
   \end{equation}
Actually, this is owed to
   $$
     \Sigma_{\phi_\alpha\in{\mathcal{S}_6(e)}}\phi_\alpha\circ(x_1x_2x_3)=\Sigma_{\phi_\alpha\in{\mathcal{S}_6(e)}}\phi_\alpha\circ(x_1^2x_2)=\Sigma_{\phi_\alpha\in{\mathcal{S}_6(e)}}\phi_\alpha\circ(x_1^3)=0.
   $$
 That's to say, there is no nontrivial $2n-$symmetric polynomial of degree three. For the case $k=2^{n}$. Note that after rotating, $\vartheta^{(k)}$ can be generated from
  $$
   e\equiv(e_1, e_2, e_3, e_4), \ \ e_1\equiv\frac{(1,1)}{\sqrt{2}}, e_2\equiv\frac{(-1,1)}{\sqrt{2}}, e_3\equiv\frac{(-1,-1)}{\sqrt{2}}, e_4\equiv\frac{(1,-1)}{\sqrt{2}}
  $$
and replacing mapping $\phi_{p,q}$ which maps $(x_1,x_2)$ to $(x_p, x_q)$ for $p<q$. Using a similar formula
   \begin{equation}\label{e5.22}
     u(x)=\frac{1}{2}\Sigma_{1\leq p<q\leq n}\Sigma_{\phi_\alpha\in{\mathcal{S}}_4(e)}\phi_{p,q}\circ\phi_\alpha\circ u_2(x_1, x_2))
   \end{equation}
for $u_2(x_1,x_2)$ stands for the component of terms of $u$ containing only variables $(x_1,x_2)$, and the fact
    \begin{equation}\label{e5.23}
    \Sigma_{\phi_\alpha\in{\mathcal{S}_4(e)}}u\big(\phi_\alpha(x_1, x_2)\big)=0,
    \end{equation}
one can infer from \eqref{e5.22}-\eqref{e5.23} that there is no $2^{n}-$symmetric polynomial of degree three. The proof of the lemma was done. $\Box$\\

Now, Theorem \ref{t5.1} is a consequence of Propositions \ref{p5.1}-\ref{p5.2} and Lemmas \ref{l5.1}-\ref{l5.5}. $\Box$

\vspace{10pt}

\section{Multiple solutions of \eqref{e1.3}}

As an application of Theorem \ref{t1.1}, we reach the following result of multiple solutions of $L_p$ dual Minkowski problem \eqref{e1.3}.

\begin{theo}\label{t6.1}
  Supposing that there is a $k-$regular polytope $\vartheta^{(k)}$ with vertices $\vartheta=(\vartheta_1,\vartheta_2,\cdots,\vartheta_k)\in{\mathbb{S}}^{n-1}$ for some $k\geq n+1$ on $n$ dimensional Euclidean space, there exists a $k-$symmetric maximal solution $h_{k}$ of the functional $F$ which is not a constant. As a result,

(1) on the planar case $n=2$, $u_k$ is of minimal period $\chi_k\equiv\frac{2\pi}{k}, \forall k\geq3$. So, if $p\leq-1$ and $q>p+(\tau+2)^2$ hold for some integer $\tau\geq1$, then \eqref{e1.3} admits at least $\tau$ solutions $\{h_k\}_{k=3}^{\tau+2}$ which are different from each other.

(2) On the higher dimensions $n\geq3$, for each $k_1\not=k_2$, $h_{k_1}$ and $h_{k_2}$ are different with each other. As a result, if $p\leq-n+1$ and
   $$
    q>p+\max_{k=k_1, \cdots, k_\tau}{\lambda_1(n,k)}
   $$
hold for different $k_i, i=1,2,\cdots,\tau$ and some integer $\tau\geq1$, \eqref{e1.3} admits at least $\tau$ different solutions $\{h_{k_i}\}_{i=1}^\tau$ up to rotations.
\end{theo}

\noindent\textbf{Proof.} On the planar case $n=2$, we note first that if the $k-$symmetric maximizer $h_k$ derived on Section 3-4 is also $l-$symmetric solution for some $l>k$, then it must be a $m_0\equiv\frac{kl}{m}-$symmetric solution for the greatest common factor of $k$ and $l$. Actually, since $k_1=\frac{k}{m}, l_1=\frac{l}{m}$ are relatively prime, for each integer $j\geq1$, there must be integers $a, b\in{\mathbb{Z}}$ such that
  $$
    j=ak_1+bl_1\Leftrightarrow \frac{2a\pi}{l}+\frac{2b\pi}{k}=\frac{2j\pi}{m_0}.
  $$
So, $h_k$ is also $m_0=\frac{kl}{m}$ symmetric. Using the fact that the Rayleigh quotient
  $$
   Q(h_k)\equiv\frac{\int^{2\pi}_0|\partial_\theta h_k|^2}{\int^{2\pi}_0h_k^2}=\frac{m^2}{l^2}Q(\overline{h_k})<\frac{\int^{2\pi}_0|\partial_\theta \overline{h_k}|^2}{\int^{2\pi}_0\overline{h_k}^2}
  $$
for
  $$
   \overline{h_k}(\theta)\equiv h_k\Big(\frac{l\theta}{m}\Big)\in\widetilde{{\mathcal{S}}}_k(\vartheta),
  $$
we arrive at the contradiction with $h_k$ is a maximizer of the functional $F$ on the family $\widetilde{{\mathcal{S}}}_k(\vartheta)$. So, $h_k$ must be different with $h_l$ for each $l>k$. Part (1) follows.

To consider the higher dimensional case, we suppose that there exist a $k$-regular polytope with vertices $\vartheta=(\vartheta_1,\cdots,\vartheta_k)\in{\mathbb{S}}^{n-1}$ satisfying
  $$
   \vartheta_1=e_1, \ \ \vartheta=\bigcup_{\phi\in{\mathcal{S}}_k(\vartheta)}\phi(e_1)
  $$
and a $\overline{k}-$regular polytope with vertices $\overline{\vartheta}=(\overline{\vartheta_1},\cdots,\overline{\vartheta_k})\in{\mathbb{S}}^{n-1}$ satisfying
  $$
   \overline{\vartheta_1}=e_1, \ \ \overline{\vartheta}=\bigcup_{\overline{\phi}\in{\mathcal{S}}_{\overline{k}}(\overline{\vartheta})}\overline{\phi}(e_1).
  $$
Then, we define the combined frame
  $$
   \vartheta^*\equiv\bigcup_{i\in{\mathbb{N}}}\widetilde{\phi}_1\circ\widetilde{\phi}_2\circ\cdots\circ\widetilde{\phi}_i(e_1), \ \ \widetilde{\phi}_i\in{\mathcal{S}}_k(\vartheta)\cup{\mathcal{S}}_{\overline{k}}(\overline{\vartheta}).
  $$
One can decompose the proof of Part (2) into the following two lemmas.

\begin{lemm}\label{l6.1}
  For $n\geq3$ and $k, \overline{k}\geq n+1$ defined above, the combined frame $\vartheta^*=\bigcup_{i}\{\vartheta^*_i\}$ spreads evenly on ${\mathbb{S}}^{n-1}$ in sense of for each $\vartheta^*_a\in\vartheta^*$, $\vartheta$ spreads evenly on the geodesic sphere $\partial B_r(\vartheta^*_a)$ for each $r>0$. Moreover, the spreading way is depending only on $r$ but not on $\vartheta^*_a$. Henceforth,

(1) if the counting number of $\vartheta^*$ equals a finite integer $m\in{\mathbb{N}}$, then $m>\max\{k,\overline{k}\}$ and there exists a $m-$regular polytope on ${\mathbb{S}}^{n-1}$.

(2) If the counting number of $\vartheta^*$ is infinite, then $\vartheta^*$ is dense on ${\mathbb{S}}^{n-1}$.
\end{lemm}

\noindent\textbf{Proof.} The evenness of $\vartheta^*$ follows from the evenness of $\vartheta$ and $\overline{\vartheta}$. Part (1) is clearly true. To show Part (2), we assume that there exists some $x\in{\mathbb{S}}^{n-1}$ such that
   \begin{equation}\label{e6.1}
     dist(x,\overline{\vartheta^*})=r_0>0.
   \end{equation}
Since the counting number of $\vartheta^*$ is infinite, there exist a sequence of distances $r_j, j\in{\mathbb{N}}$ and a sequence of $\vartheta^*_j\in\vartheta^*$, such that
   \begin{equation}\label{e6.2}
    dist(\vartheta^*_a,\vartheta^*_j)=r_j>0, \ \ \lim_{j\to\infty}r_j=0
   \end{equation}
holds for arbitrarily chosen $\vartheta^*_a\in\vartheta^*$. Then, for $\varepsilon>0$ small, there exists some $\vartheta^*_\varepsilon\in\vartheta^*$ such that
   \begin{equation}\label{e6.3}
     dist(x,\vartheta^*_\varepsilon)\leq r_0+\varepsilon, \ \ dist(x,\vartheta^*)=r_0.
   \end{equation}
However, it is inferred from \eqref{e6.2} that the set
   $$
    {\mathcal{C}}_{\varepsilon,j}\equiv\vartheta^*\cap\partial B_{r_j}(\vartheta^*_\varepsilon)
   $$
is nonempty and spreads evenly on the geodesic sphere $\partial B_{r_j}(\vartheta^*_\varepsilon)$. This contradicts with \eqref{e6.3} for $j$ large and then $\varepsilon$ small, since the empty part of ${\mathcal{C}}_{\varepsilon,j}$ contains almost half geodesic sphere $\partial B_{r_j}(\vartheta^*_\varepsilon)$. The proof was done. $\Box$

\begin{lemm}\label{l6.2}
   For dimension $n\geq3$, we assume that $h_{k}$ is a $k-$symmetric maximizer of $F$ on the family ${\mathcal{C}}_k$. Then $h_{k}$ can not be a $\overline{k}-$symmetric function for another $\overline{k}\not=k$.
\end{lemm}

\noindent\textbf{Proof.} As in the Lemma \ref{l6.1}, if case (1) occurs, then $\vartheta^*$ forms a $m-$regular polytope containing $\vartheta$ as a $k-$regular sub-polytope for $m$ being divisible by $k$. However, this is impossible on dimensions $n\geq3$ by Proposition \ref{p2.1}. Supposing that $h_{\vartheta^{(k)}}$ is also a $\overline{k}-$symmetric function for another polytope $\vartheta^{(\overline{k})}$, if case (2) occurs, $h_{\vartheta^{(k)}}$ must be a constant function on ${\mathbb{S}}^{n-1}$. This contradicts with Theorem \ref{t1.1} and thus gives the proof of Part (2) of Theorem \ref{t6.1}. $\Box$\\

Now, Part (2) of Theorem \ref{t6.1} was a direct consequence of Lemma \ref{l6.2}. The proof of Theorem \ref{t6.1} was done. $\Box$

\vspace{10pt}

\section{Chou-Wang type identity for \eqref{e1.2}}

From now on, we turn to consider solvability problem of \eqref{e1.2}. Given an arbitrarily $n\times n$ matrix $B$, it generates a projective vector field
     \begin{equation}\label{e7.1}
       \xi(x)\equiv Bx-(x^TBx)x, \ \ x\in{\mathbb{S}}^{n-1}.
     \end{equation}
Writing $\xi=\xi^ke_k$ for $\xi^k=\langle \xi,e_k\rangle=e_k^TBx$, we note first that
   \begin{eqnarray}\nonumber
     &&\xi^k_{,i}=x^T_{,k}Bx_{,i}-x^TBx\delta_{ki},\\
     &&\xi^k_{,ij}=-\xi^k\delta_{ij}-\xi^j\delta_{ki}-x^TBx_{,i}\delta_{jk}-x^TBx_{,j}\delta_{ki}.
   \end{eqnarray}
Now, let us now prove a key lemma would be useful later.

\begin{lemm}\label{l7.1}
  For any $C^3-$function $h$ on ${\mathbb{S}}^{n-1}$, $\beta\in{\mathbb{R}}$ and any tangent vector field $\xi$ of the form \eqref{e7.1}, there holds
    \begin{eqnarray}\nonumber
     &&n\int_{{\mathbb{S}}^{n-1}}|Z|^\beta\det(A)\nabla_\xi h=\int_{{\mathbb{S}}^{n-1}}|Z|^\beta h\det(A)div\xi-\beta\int_{{\mathbb{S}}^{n-1}}|Z|^{\beta-2}\det(A)|\nabla h|^2\xi^kh_k\\
     &&\ \ \ \ \ \ \ \ \ +\beta\int_{{\mathbb{S}}^{n-1}}|Z|^{\beta-2}\det(A)hh_jh_k\xi^k_{,j}+\beta\int_{{\mathbb{S}}^{n-1}}|Z|^{\beta-2}h^2\det(A)h_k(x^TBx_{,k}).
    \end{eqnarray}
\end{lemm}

\noindent A special case of Lemma \ref{l7.1} for $\beta=0$ has been shown in \cite{lu2013rotationally}.

\noindent\textbf{Proof.} Taking any function $\eta$,
  \begin{eqnarray*}
   (n-1)\int_{{\mathbb{S}}^{n-1}}\eta\det(A)\nabla_\xi h&=&\int_{{\mathbb{S}}^{n-1}}\eta U^{ij}(h_{ij}+h\delta_{ij})\xi^kh_k=\int_{{\mathbb{S}}^{n-1}}U^{ij}h(\eta\xi^kh_k)_{,ij}+\eta U^{ij}h\delta_{ij}\xi^kh_k\\
   &=&\int_{{\mathbb{S}}^{n-1}}\eta hU^{ij}h_k(\xi^k_{,ij}+\xi^k\delta_{ij})+2\eta hU^{ij}\xi^k_{i,i}h_{,kj}+\eta h U^{ij}\xi^k(A_{ki,j}-h_j\delta_{ki})+R_\eta^{(1)}
  \end{eqnarray*}
holds for
  \begin{equation}
   R^{(1)}_\eta\equiv2\int_{{\mathbb{S}}^{n-1}}U^{ij}h\eta_i(\xi^k_{,j}h_k+\xi^kh_{kj})+\int_{{\mathbb{S}}^{n-1}}U^{ij}h\eta_{ij}\xi^kh_k.
  \end{equation}
Using the symmetry of $A_{ki,j}=A_{ij,k}$, one can proceed further that
   \begin{eqnarray*}
    (n-1)\int_{{\mathbb{S}}^{n-1}}\eta\det(A)\nabla_\xi h&=&-\int_{{\mathbb{S}}^{n-1}}\eta hU^{ij}h_k(\xi^j\delta_{ki}+x^TBx_{,i}\delta_{jk}+x^TBx_{,j}\delta_{ki})\\
    &&+\int_{{\mathbb{S}}^{n-1}}2\eta hU^{ij}\xi^k_{,i}h_{kj}+\eta hU^{ij}\xi^kA_{ij,k}-\eta hU^{ij}h_j\xi^i+R^{(1)}_\eta\\
    &=&J_1-\int_{{\mathbb{S}}^{n-1}}\det(A)(\eta h\xi^k)_{,k}+R^{(1)}_\eta,
   \end{eqnarray*}
where
   \begin{eqnarray*}
    J_1&\equiv&\int_{{\mathbb{S}}^{n-1}}-2\eta hU^{ij}h_j\xi^i-2\eta hU^{ij}h_jx^TBx_{,i}+2\eta hU^{ij}\xi^k_{,i}h_{,kj}\\
    &=&\int_{{\mathbb{S}}^{n-1}}h^2U^{ij}(\eta\xi^i+\eta x^TBx_{,i})_{,j}+2\eta hU^{ij}\xi^k_{,i}h_{,kj}\\
    &=&2\int_{{\mathbb{S}}^{n-1}}\eta h\det(A)div\xi+\int_{{\mathbb{S}}^{n-1}}h^2U^{ij}\eta_j(\xi^i+x^TBx_{,i}).
   \end{eqnarray*}
Henceforth,
  \begin{eqnarray}\nonumber
   n\int_{{\mathbb{S}}^{n-1}}\eta\det(A)\nabla_\xi h&=&\int_{{\mathbb{S}}^{n-1}}\eta h\det(A)div\xi+R^{(2)}_\eta
  \end{eqnarray}
holds for
  \begin{eqnarray}\nonumber\label{e7.5}
   R^{(2)}_\eta&\equiv&    R^{(1)}_\eta-\int_{{\mathbb{S}}^{n-1}}h\det(A)\eta_k\xi^k+\int_{{\mathbb{S}}^{n-1}}h^2U^{ij}\eta_j(\xi^i+x^TBx_{,i})\\
\nonumber   &=&2\int_{{\mathbb{S}}^{n-1}}U^{ij}h\eta_i(\xi^k_{,j}h_k+\xi^kh_{kj})+\int_{{\mathbb{S}}^{n-1}}U^{ij}h\eta_{ij}\xi^kh_k\\
   &&-\int_{{\mathbb{S}}^{n-1}}h\det(A)\eta_k\xi^k+\int_{{\mathbb{S}}^{n-1}}h^2U^{ij}\eta_j(\xi^i+x^TBx_{,i})\\
\nonumber   &=&2\int_{{\mathbb{S}}^{n-1}}U^{ij}hh_k\eta_i\xi^k_{,j}+\int_{{\mathbb{S}}^{n-1}}h\det(A)\xi^k\eta_k-\int_{{\mathbb{S}}^{n-1}}h^2U^{ij}\eta_i\xi^j\\
 \nonumber  &&+\int_{{\mathbb{S}}^{n-1}}hU^{ij}\eta_{ij}\xi^kh_k+\int_{{\mathbb{S}}^{n-1}}h^2U^{ij}\eta_j(x^TBx_{,i}).
  \end{eqnarray}
Now, we take $\eta=|Z|^\beta$ for $Z=h(x)x+\nabla h$ and use the relation
  \begin{eqnarray}\nonumber
   \eta_i&=&\beta|Z|^{\beta-2}\langle Z_{,i},Z\rangle=\beta|Z|^{\beta-2}A_{ik}h_k\\
   \eta_{ij}&=&\beta(\beta-2)|Z|^{\beta-4}A_{ik}h_kA_{jl}h_l+\beta|Z|^{\beta-2}A_{ik,j}h_k+\beta|Z|^{\beta-2}A_{ik}h_{kj}.
  \end{eqnarray}
Substituting into \eqref{e7.5}, we rewrite
  \begin{eqnarray*}
   R^{(2)}_\eta&=&2\beta\int_{{\mathbb{S}}^{n-1}}|Z|^{\beta-2}\det(A)hh_jh_k\xi^k_{,j}+\beta\int_{{\mathbb{S}}^{n-1}}|Z|^{\beta-2}h\det(A)A_{kl}h_l\xi^k\\
   &&-\beta\int_{{\mathbb{S}}^{n-1}}|Z|^{\beta-2}h^2\det(A)\xi^kh_k+\beta(\beta-2)\int_{{\mathbb{S}}^{n-1}}|Z|^{\beta-4}h\xi^kh_k\det(A)A_{ij}h_ih_j\\
   &&+\beta\int_{{\mathbb{S}}^{n-1}}|Z|^{\beta-2}h[\det(A)]_{,l}h_l\xi^kh_k+\beta\int_{{\mathbb{S}}^{n-1}}|Z|^{\beta-2}h\det(A)\triangle h\xi^kh_k\\
   &&+\beta\int_{{\mathbb{S}}^{n-1}}|Z|^{\beta-2}h^2\det(A)h_k(x^TBx_{,k}).
  \end{eqnarray*}
By divergence theorem and $|Z|_{,i}=|Z|^{-1}A_{ik}h_k$, one can simplify further that
  \begin{eqnarray*}
   R^{(2)}_\eta&=&\beta\int_{{\mathbb{S}}^{n-1}}|Z|^{\beta-2}\det(A)hh_jh_k\xi^k_{,j}+\beta\int_{{\mathbb{S}}^{n-1}}|Z|^{\beta-2}h^2\det(A)h_k(x^TBx_{,k})\\
   &&-\beta\int_{{\mathbb{S}}^{n-1}}|Z|^{\beta-2}\det(A)|\nabla h|^2\xi^kh_k.
  \end{eqnarray*}
The proof was done. $\Box$\\

As an application of Lemma \ref{l7.1}, one can deduce the following Chou-Wang type identity ensuring the solvability of \eqref{e1.2}. The identity of this form was firstly obtained by Chou and Wang in \cite{chou2006lp} for $q=n$ and $p=-n$. The version we presented here includes all exponents $pq\not=0$.

\begin{theo}\label{t7.1}
  For any $C^2-$solution $h$ of \eqref{e1.2} for some Lipschitz nonnegative function $f$ on ${\mathbb{S}}^{n-1}$, and any tangent vector field $\xi$ of the form \eqref{e7.1}, there holds
  \begin{equation}\label{e7.7}
    \int_{{\mathbb{S}}^{n-1}}\Big(\nabla_\xi f+\frac{p+q}{q}fdiv\xi\Big)h^p=-\frac{p(q-n)}{q}\int_{{\mathbb{S}}^{n-1}}fh^p|Z|^{-2}(h\nabla_\xi h+\xi^k_{,j}h_jh_k).
  \end{equation}
\end{theo}

\noindent\textbf{Proof.} Taking $\beta=q-n$ in Lemma \ref{l7.1} and noting that $\xi^k=x^TBx_{,k}$, by the equation \eqref{e1.2}, there holds
  \begin{eqnarray*}
   q\int_{{\mathbb{S}}^{n-1}}fh^{p-1}\nabla_\xi h&=&\int_{{\mathbb{S}}^{n-1}}fh^pdiv\xi+2(q-n)\int_{{\mathbb{S}}^{n-1}}fh^{p+1}|Z|^{-2}\nabla_\xi h\\
   &&+(q-n)\int_{{\mathbb{S}}^{n-1}}fh^p|Z|^{-2}\xi^k_{,j}h_jh_k.
  \end{eqnarray*}
\eqref{e7.7} follows from integration by parts. $\Box$\\

As a corollary of Theorem \ref{t7.1}, we obtain the following solvability condition.

\begin{coro}\label{c7.1}
  Letting $h$ be a $C^2({\mathbb{S}}^{n-1})-$solution of \eqref{e1.2} for Lipschitz nonnegative function $f$ on ${\mathbb{S}}^{n-1}$, we have
   \begin{equation}\label{e7.8}
      \int_{{\mathbb{S}}^{n-1}}\Bigg(\nabla_\xi f+\frac{p+q}{q}fdiv\xi+\frac{p(q-n)}{q}f\chi(\xi)\Bigg)h^p=0
   \end{equation}
holds for
    $$
     \chi(\xi)\equiv\frac{h\nabla_\xi h+B(\nabla h,\nabla h)-(trB-div\xi)/n\times|\nabla h|^2}{|Z|^2},
    $$
where $\lambda_1\geq\lambda_2\geq\cdots\geq\lambda_{n}$ are eigenvalues of $B$ and
   \begin{equation}\label{e7.9}
   \begin{cases}
    \chi(\xi)\leq|\xi|+\sup_iB(e_i,e_i)-\frac{trB-div\xi}{n}\leq2\lambda_1(B)-\lambda_{n}(B),\\[5pt]
    \chi(\xi)\geq-|\xi|+\inf_iB(e_i,e_i)-\frac{trB-div\xi}{n}\geq-2\lambda_1(B)+\lambda_{n}(B).
   \end{cases}
  \end{equation}
\end{coro}

\noindent\textbf{Proof.} Noting that
 $$
  \xi^k_{,j}h_jh_k=e_j^TBe_kh_jh_k-(x^TBx)|\nabla h|^2
 $$
and
 $$
   div\xi=e_i^TBe_i-(n-1)x^TBx=trB-nx^TBx,
 $$
\eqref{e7.9} follows from definitions of $B, \xi$ and the observation
   $$
    \lambda_n(B)\leq\frac{trB-div\xi}{n}\leq\lambda_1(B).
   $$
 $\Box$

\vspace{10pt}

\section{Nonexistence result for $p\leq-q, q\geq n$}

A first application of Theorem \ref{t7.1} is the following nonexistence result for $p=-q$.

\begin{theo}\label{t8.1}
  Supposing that $n\geq2$ and $p=-q, q\geq n$, then there exists a smooth function $f$ which is positive outside equator or outside two polar, such that \eqref{e1.2} has no solution.
\end{theo}

\noindent\textbf{Proof.} Taking $B=diag(0,\cdots,0,1)$, it generates a projective vector field
   $$
    \xi=Bx-(x^TBx)x=\left(
      \begin{array}{c}
       -x_{n}^2x'\\
       x_{n}(1-x_{n}^2)
      \end{array}
      \right)
   $$
for $x=(x',x_{n})\in{\mathbb{S}}^{n-1}$ satisfying
   \begin{equation}\label{e8.1}
     \begin{cases}
        div\xi=1-nx_{n}^2,\ \ |\xi|=x_{n}^2(1-x_{n}^2),\\
       \chi(\xi)\leq x_{n}^2(1-x_{n}^2)+1-x_{n}^2=1-x_{n}^4\\
       \chi(\xi)\geq-x_{n}^2(1-x_{n}^2)-x_{n}^2=-x_{n}^2(2-x_{n}^2).
    \end{cases}
   \end{equation}
When $q\geq n$, one can construct a smooth function
  $$
   f(x)=x_{n}^{2\kappa}, \ \ \forall x_{n}\in[-1,1]
  $$
 for some $\kappa=[q-n]_*$ and $[z]_*$ stands for the least integer greater than $z$, which is positive outside equator. Noting that
  $$
    \nabla_\xi f=x_{n}(1-x_{n}^2)\frac{df}{dx_{n}},
  $$
we have
  \begin{eqnarray*}
   &&\nabla_\xi f-(q-n)f\chi(\xi)\geq x_{n}(1-x_{n}^2)\frac{df}{dx_{n}}-(q-n)(1-x_{n}^4)f\\
   &&\ \ \ \ \geq(1-x_{n}^2)\Big\{2\kappa-(q-n)(1+x_{n}^2)\Big\}f\geq0
  \end{eqnarray*}
holds everywhere on ${\mathbb{S}}^{n-1}$. So, we conclude a contradiction from \eqref{e7.8} as long as a $C^2$-solution $h$ exists. To construct a smooth function $f$ which is positive outside two polar, we construct a smooth function
  $$
   f(x)=(1-x_{n}^2)^{\kappa}, \ \ \kappa\equiv[q-n]_*, \ \ \forall x_{n}\in[-1,1].
  $$
Using the lower bound of $\chi(\xi)$, we have
  \begin{eqnarray*}
   &&\nabla_\xi f-(q-n)f\chi(\xi)\leq x_{n}(1-x_{n}^2)\frac{df}{dx_{n}}+(q-n)x_{n}^2(2-x_{n}^2)f\\
   &&\ \ \ \ \leq[-2\kappa x_{n}^2+2(q-n)x_{n}^2]f\leq0, \ \ \forall x\in{\mathbb{S}}^{n-1}.
  \end{eqnarray*}
A contradiction occurs again from \eqref{e7.8}. $\Box$\\

When $p<-q$, there is one more term $\frac{p+q}{q}fdiv\xi$ appears on the identity \eqref{e7.8}. After handling this term carefully, we have also the following nonexistence result.

\begin{theo}\label{t8.2}
 Supposing that $n\geq2$ and $p<-q, q\geq n$, then there exists a smooth function $f$ which is positive outside equator or two polar, such that \eqref{e1.2} has no solution.
\end{theo}

\noindent\textbf{Proof.} Choosing $B=diag(0,\cdots,0,1)$ as in proof of Theorem \ref{t8.1}, and using the relations \eqref{e8.1}, there exists a large integer $\gamma_{p,q}\geq1$ such that
  \begin{eqnarray}\nonumber\label{e8.2}
   &&K_f\equiv\nabla_\xi f+\frac{p+q}{q}fdiv\xi+\frac{p(q-n)}{q}f\chi(\xi)\\
   && \ \ \ \ \ \ \ \leq x_{n}(1-x_{n}^2)\frac{df}{dx_{n}}+\gamma_{p,q}\Big[\Big(x_{n}^2-\frac{1}{n}\Big)_++x_{n}^2\Big]f.
  \end{eqnarray}
Taking another large integer $\beta\gg\gamma_{p,q}$, we construct
  $$
   f(x)=e^{-\beta x_{n}^2}(1-x_{n}^2)^\beta, \ \ \forall x\in{\mathbb{S}}^{n-1}
  $$
to be a smooth nonnegative function which is positive outside two polar. Then, the inequality
  \begin{equation}
   K_f\leq f\Big\{-2\beta x_{n}^2(1-x_{n}^2)-2\beta x_{n}^2+\gamma_{p,q}\Big[\Big(x_{n}^2-\frac{1}{n}\Big)_++x_{n}^2\Big]\Big\}\leq0
  \end{equation}
holds for all $x\in{\mathbb{S}}^{n-1}$, as long as $\beta$ is sufficiently large. We thus reach a contradiction. If one wants to construct nonnegative function $f$ which is positive outside equator, we replace \eqref{e8.2} by
  \begin{eqnarray}\nonumber
   &&K_f\equiv\nabla_\xi f+\frac{p+q}{q}fdiv\xi+\frac{p(q-n)}{q}f\chi(\xi)\\
   && \ \ \ \ \ \ \ \geq x_{n}(1-x_{n}^2)\frac{df}{dx_{n}}-\gamma_{p,q}\Big[\Big(\frac{1}{n}-x_{n}^2\Big)_++(1-x_{n}^2)\Big]f.
  \end{eqnarray}
Taking a large integer $\beta\gg\gamma_{p,q}$, we construct
  $$
   f(x)=x_{n}^{2\beta}, \ \ \forall x\in{\mathbb{S}}^{n-1}
  $$
to be a smooth nonnegative function which is positive outside equator. Then, the inequality
  \begin{equation}
   K_f\geq f\Big\{2\beta-\gamma_{p,q}\Big[\Big(\frac{1}{n}-x_{n}^2\Big)_++(1-x_{n}^2)\Big]\Big\}\geq0
  \end{equation}
holds for all $x\in{\mathbb{S}}^{n-1}$, as long as $\beta$ is sufficiently large. We still have a contradiction. $\Box$


\end{document}